\documentclass[12pt,a4paper,draft]{article}
\usepackage[T2A]{fontenc}
\usepackage{latexsym,amsfonts,amssymb,amsmath,longtable,amsthm,mathrsfs}
\renewcommand{\le}{\leqslant}
\renewcommand{\ge}{\geqslant}

\newtheorem{Lemma}{{\bfseries Lemma}}
\newtheorem{Cor}[Lemma]{{\bfseries Corollary}}
\newtheorem{Theo}[Lemma]{{\bfseries Theorem}}
\newtheorem{Hyp}[Lemma]{{\bfseries Hypothesis}}

\DeclareMathOperator{\GL}{GL} \DeclareMathOperator{\Aut}{Aut}

\title{\vspace{-1cm} \hfill{\normalsize MSC2010 20D20, 20E34}{
\fontfamily{cmr} \fontseries{bx} \selectfont \\ \vspace{1cm} A conjugacy criterion for Hall subgroups in finite groups\footnote{The authors
were supported by  Russian Foundation for Basic Research (Grants 08-01-00322 and 10-01-00391),  the Russian Federal Agency for Education
(Grant 2.1.1.419), and
Federal Target Grant (Contract No. 02.740.11.0429). The first author gratefully acknowledges the support from Deligne 2004 Balzan prize
in mathematics, and the Lavrent'ev Young Scientists Competition (No 43 on
04.02.2010). }}
\author{\bf  E.P. Vdovin,
D.O. Revin}}

\begin{document}
%%%%%%%%%%%%%%%%%%%%%%%%%%%%%%%%%%%%%%%%%%%%%%%%%%%%%%%%%%%%%%%%%%%%%%%%%%%%%%%%%%%

%\renewcommand{\qed}{}
%\renewcommand{\proofname}{{ {\rm ДОКАЗАТЕЛЬСТВО}}}

%\thispagestyle{empty}
%\vspace{-2cm}

\maketitle
\pagenumbering{arabic}
\begin{abstract}

A finite group $G$ is said to satisfy $C_\pi$ for a set of primes $\pi$, if $G$ possesses exactly one
class of conjugate $\pi$-Hall subgroups. In the paper we obtain a criterion for a finite group $G$ to satisfy
$C_\pi$ in terms of a  normal series of the group.
%\medskip

\end{abstract}

%%%%%%%%%%%%%%%%%%%%%%%%%%%%%%%%%%%%%%%%%%%%%%%%%%%%%%%%%%%%%%%%%%%%%%%

\section*{Introduction}

%%%%%%%%%%%%%%%%%%%%%%%%%%%%%%%%%%%%%%%%%%%%%%%%%%%%%%%%%%%%%%%%%%%%%%%%
Let  $\pi$ be a set of primes. We denote by  $\pi'$ the set of all primes not in $\pi$, by $\pi(n)$ the set of
prime divisors of a positive integer $n$, while for a finite group  $G$ by $\pi(G)$ we denote $\pi(|G|)$. A
positive integer $n$ with $\pi(n)\subseteq\pi$ is called a {\it $\pi$-number}, while a group $G$ with
$\pi(G)\subseteq \pi$ is called a {\it $\pi$-group}. A subgroup $H$ of $G$ is called a {\it $\pi$-Hall subgroup},
if $\pi(H)\subseteq\pi$ and $\pi(|G:H|)\subseteq \pi'$.  According to [1]  we  say that  $G$ {\it satisfies
$E_\pi$} (or briefly $G\in E_\pi$), if $G$ possesses a $\pi$-Hall subgroup.  Moreover, if every two
$\pi$-Hall subgroups are conjugate, then we say that $G$ {\it satisfies $C_\pi$} ($G\in C_\pi$). 
Further if  each $\pi$-subgroup of $G$ lies in a $\pi$-Hall subgroup then we say that  $G$ {\it
satisfies $D_\pi$} ($G\in D_\pi$). A group satisfying $E_\pi$ (respectively, $C_\pi$, $D_\pi$)  we also call an
$E_\pi$- (respectively,  $C_\pi$-, $D_\pi$-) {\it group}.

Let $A,B$, and $H$ be subgroups of $G$ such that $B\unlhd A$. We denote $N_H(A)\cap N_H(B)$ by $N_H(A/B)$. Then every
element $x\in N_H(A/B)$ induces an automorphism of $A/B$ acting by $Ba\mapsto B x^{-1}ax$.
Thus the homomorphism $N_H(A/B)\rightarrow \operatorname{Aut}(A/B)$ is defined. The image of the homomorphism is
denoted by $\operatorname{Aut}_H(A/B)$ and is called the {\it group of $H$-induced automorphisms on $A/B$}, while
the kernel is denoted by~${C_H(A/B)}$. If $B=1$, then $\operatorname{Aut}_H(A/B)$ is denoted
by~$\operatorname{Aut}_H(A)$.

Assume that $\pi$ is fixed. It is proved that the class of  $D_\pi$-groups is closed under homomorphic
images, normal subgroups (mod CFSG\footnote{(mod CFSG) in the paper means that the result is proven
modulo the classification of finite simple groups}), [2, Theorem~7.7] or [3, Corollary~1.3]), and extensions (mod
CFSG, [2, Theorem~7.7]). Thus, a finite group $G$ satisfies $D_\pi$ if and only if each composition factor $S$ of
$G$ satisfies $D_\pi$. The class of $E_\pi$-groups is also known to be closed under normal subgroups and
homomorphic images (see Lemma~4(1)), but not closed under extensions in general (see [4, Ch.~V,
Example~2]). In  [5, Theorem~3.5] and [6, Corollary~6] it is proven that, if $1=G_0<G_1<\ldots<G_n=G$ is a
composition series of $G$ that is a refinement of a chief series, then   $G$ satisfies $E_\pi$ if and only if
$\operatorname{Aut}_G(G_i/G_{i-1})$ satisfies $E_\pi$ for each~${i=1,\ldots,n}$.

The class of $C_\pi$-groups is closed under extensions (see Lemma~5) but not closed under normal subgroups in general (see the example
below). In the present paper, by using the classification of finite simple groups, we show that the class of $C_\pi$-groups is closed under
homomorphic images (mod CFSG, see Lemma~9), and give a criterion for a finite group to satisfy $C_\pi$ in terms of a normal series of the
group. The main result is the following:

\begin{Theo}\label{MainTheorem} {\em (mod CFSG)}
Let $\pi$ be a set of primes, let $H$ be a $\pi$-Hall subgroup, and $A$ be a normal subgroup of a $C_\pi$-group $G$. Then~${HA\in C_\pi}$.
\end{Theo}

\begin{Cor}\label{FirstCorollary} {\em (a conjugacy criterion for Hall subgroups, mod CFSG).}
Let $\pi$ be a set of
primes and let $A$ be a normal subgroup of $G$. Then $G\in C_\pi$ if and only if $G/A\in C_\pi$ and, for a $\pi$-Hall
\footnote{Since $G/A\in C_\pi $ the phrase ``for a $\pi$-Hall subgroup $K/A$ of $G/A$'' in the
statement can be interpreted as ``for every'' and ``for some'', and both of them are correct.} subgroup $K/A$ of
$G/A$ its complete preimage  $K$ satisfies $C_\pi$. In particular, if  $\vert G:A\vert$ is a $\pi'$-number, then
$G\in C_\pi$ if and only if~${A\in C_\pi}$.
\end{Cor}

Using the corollary, at the end of the paper we give an algorithm that reduces the problem whether a finite groups satisfies $C_\pi$
to the same problem in some almost simple groups. In view of Theorem~1 notice that we do not know any counterexample to the following
hypothesis.

\begin{Hyp}
Let $\pi$ be a set of primes and let $A$ be a {\rm (}not necessary normal{\rm )} subgroup of a finite $C_\pi$-group $G$ containing a
$\pi$-Hall subgroup of $G$. Then~${A\in C_\pi}$.
\end{Hyp}

In the hypothesis the condition that $A$ includes a  $\pi$-Hall subgroup of $G$ cannot be weaken by the condition that
the index of $A$ is a  $\pi'$-number. Indeed, consider $B_3(q)\simeq \mathrm{P}\Omega_7(q)$, where $q-1$ is divisible by $12$ and is not
divisible by $8$ and $9$. In view of [2, Lemma~6.2],
$\mathrm{P}\Omega_7(q)$ is a $C_{\{2,3\}}$-group and its
$\{2,3\}$-Hall subgroup is included in a monomial subgroup. On the other hand, $\Omega_7(2)$ is known to be isomorphic to a subgroup of
$\mathrm{P}\Omega_7(q)$ and, under the above stated conditions on $q$, its index is not divisible by $2$ and $3$. However
$\Omega_7(2)$ does not possess a $\{2,3\}$-Hall subgroup, i.e., it is not even an $E_{\{2,3\}}$-group.

\section{Notations and preliminary results}

By $\pi$  we always denote a set of primes, and the term ``group'' always means a finite group.

The lemmas below  are known and their proof do not refer to the classification of finite simple groups.

\begin{Lemma} {\rm [4, Ch.~IV, (5.11), Ch.~V, Theorem~3.7]}
Let $A$ be a normal subgroup of $G$. Then the following holds.

{\rm (1)} If $H$ is a $\pi$-Hall subgroup of $G$ then
$H \cap A$ is a $\pi$-Hall subgroup of $A$, while $HA/A$ is a
$\pi$-Hall subgroup of~${G/A}$.

{\rm (2)} If all factors of a subnormal series of $G$ are either $\pi$- or
$\pi'$- groups then~${G\in D_{\pi}}$.
\end{Lemma}

Notice that (2) of Lemma~4 follows from the  famous Chunikhin's theorem on
$\pi$-solvable groups and the Feit-Thompson Odd Order Theorem.

\begin{Lemma} {\rm (Chunikhin; see also [1, Theorems~C1
and~C2] or [4, Ch. V, (3.12)])}
Let $A$ be a subnormal subgroup of $G$. If both $A$ and $G/A$ satisfy $C_{\pi}$, then~${G\in C_{\pi}}$.
\end{Lemma}

\begin{Lemma}{\rm [3, Lemma~2.1(e)]} Let $A$ be a normal subgroup of $G$  such that $G/A$ is a $\pi$-group and
let $M$ be a $\pi$-Hall subgroup of $A$. Then a $\pi$-Hall subgroup $H$ of $G$ with $H\cap A=M$ exists if and
only if $G$ acting by conjugation leaves $\{M^a\mid a\in A\}$ invariant.
\end{Lemma}

{\scshape Example.} Suppose that $\pi=\{2,3\}$. Let $G=\operatorname{GL}_5(2)=\operatorname{SL}_5(2)$ be a group of order
$99999360=2^{10}\cdot 3^2\cdot 5\cdot 7 \cdot 31$. Assume then that $\iota : x\in G\mapsto (x^t)^{-1}$ and
$\widehat{G}=G\leftthreetimes \langle\iota\rangle$ is a natural semidirect product. By [7, Theorem~1.2],  $G$
possesses $\pi$-Hall subgroups and each $\pi$-Hall subgroup is the stabilizer of a series of subspaces
$V=V_0<V_1<V_2<V_3=V$, where $V$ is a natural module of $G$ and $\dim V_k/V_{k-1}\in \{1,2\}$ for every
$k=1,2,3$. Hence $G$ possesses exactly three classes of conjugate $\pi$-Hall subgroups with representatives
$$
%\xymatrix
H_1=\left(
\begin{array}{c@{}c@{}c}
\fbox{$\begin{array}{c}
\\
\!
\GL_2(2)
\!
\\
\\
\end{array}$}& &*\\
 &\fbox{1}& \\
0& &
\fbox{$\begin{array}{c}
\\
\!
\GL_2(2)
\!
\\
\\
\end{array}$}
\end{array}
\right),
$$
$$H_2=
{\left(
\begin{array}{c@{}c@{}c}
\fbox{1}& &*\\
 &\fbox{$\begin{array}{c}
\\
\!
\GL_2(2)
\!
\\
\\
\end{array}$}& \\
0& &
\fbox{$\begin{array}{c}
\\
\!
\GL_2(2)
\!
\\
\\
\end{array}$}
\end{array}
\right),
%\stackrel{\iota}{\longleftrightarrow}
\text{ and }
H_3=\left(
\begin{array}{c@{}c@{}c}
\fbox{$\begin{array}{c}
\\
\!
\GL_2(2)
\!
\\
\\
\end{array}$}& &*\\
 &\fbox{$\begin{array}{c}
\\
\!
\GL_2(2)
\!
\\
\\
\end{array}$}& \\
0& &
\fbox{1}
\end{array}
\right).}
$$
 Notice that
$N_G(H_k)=H_k$, $k=1,2,3$, since $H_k$ is parabolic. By Lemma~4(1), for every $\pi$-Hall subgroup $H$ of
$\widehat{G}$ the subgroup $H\cap G$ is conjugate to one of subgroups  $H_1, H_2$, and $H_3$. The class containing
$H_1$ is invariant under $\iota$. Hence by Lemma~6 there exists a $\pi$-Hall subgroup $H$ of $\widehat{G}$ with
$H\cap G=H_1$. Moreover $H=N_{\widehat{G}}(H_1)$. The classes containing $H_2$ and $H_3$ are permuted by~$\iota$.
So from Lemmas~4(1) and~6 it follows that these subgroups do not lie in $\pi$-Hall subgroups of
$\widehat{G}$. Thus $\widehat{G}$ has exactly one class of conjugate $\pi$-Hall subgroups, and so satisfies $C_\pi$,
in contrast to its normal subgroup $G$.

\begin{Lemma} Let $A$ be a normal subgroup and let $H$ be a $\pi$-Hall subgroup of a $C_\pi$-group $G$. Then 
$N_G(HA)$ and $N_G(H\cap A)$ satisfy~$C_\pi$. \end{Lemma}

\begin{proof} Note that both  $N_G(HA)$ and $N_G(H\cap A)$ include $H$, and so they satisfy $E_\pi$. Let $K$ be a
$\pi$-Hall subgroup of $N_G(HA)$. Since  $HA\unlhd N_G(HA)$ and $|N_G(HA):HA|$ is a $\pi'$-number, we have $K\leq
HA$ and $KA=HA$. If $x\in G$ is chosen so that $K=H^x$ then $(HA)^x=H^xA=KA=HA$ and so $x\in N_G(HA)$. Therefore
\hbox{$N_G(HA)\in C_\pi$.}

Suppose that $K$ is a $\pi$-Hall subgroup of $N_G(H\cap A)$. Then \hbox{$K(H\cap A)=K$} and $K\cap A=H\cap A$. If
$x\in G$ is chosen so that $K=H^x$, then $(H\cap A)^x=H^x\cap A=K\cap A=H\cap A$ and so $x\in N_G(H\cap A)$.
Therefore $N_G(H\cap A)\in C_\pi$. \end{proof}

\begin{Lemma}{\rm ([6, Corollary~9] mod CFSG)} Each $\pi$-Hall of a homomorphic image of an $E_\pi$-gro\-up
$G$ is the image of a $\pi$-Hall subgroup of~$G$. \end{Lemma}

By Lemma 8 it is immediate that $C_\pi$ is preserved under homomorphisms.

\begin{Lemma}{\rm(mod CFSG)} Let $A$ be a normal subgroup of a $C_\pi$-group $G$. Then  $G/A\in C_\pi$.
\end{Lemma}

\begin{proof} Put $G/A=\overline{G}$. Since all $\pi$-Hall subgroups of $G$ are conjugate, it is enough to show
that for every $\pi$-Hall subgroup $\overline{K}$ of $\overline{G}$ there exists a $\pi$-Hall subgroup $U$ of $G$
such that $UA/A=\overline{K}$. The existence of $U$ follows from Lemma~8. \end{proof}

If $G$ -is a group then a {\it $G$-class of $\pi$-Hall subgroups} is a class of conjugate $\pi$-Hall subgroups
of $G$. Let $A$ be a subnormal subgroup of a $E_\pi$-group $G$. A subgroup  $H\cap A$ of $A$, where $H$ is a
$\pi$-Hall subgroup of $G$ is called a {\it $G$-induced $\pi$-Hall subgroup} of $A$. Thus the set $\{(H\cap
A)^a\mid a\in A\}$, where $H$ is a $\pi$-Hall subgroup of $G$ is called an {\it $A$-class of $G$-induced
$\pi$-Hall subgroups}. We denote the number of all $A$-classes of $G$-induced $\pi$-Hall subgroups by
$k_\pi^G(A)$. Let $k_\pi(G)=k_\pi^G(G)$ be the number of classes of $\pi$-Hall subgroups of~$G$.
Clearly~${k_\pi^G(A)\le k_\pi (A)}$. 

Recall that a finite group $G$ is called {\it almost simple}, if $G$ possesses a unique minimal normal subgroup
$S$ and $S$ is a nonabelian finite group (equivalently, up to isomorphism, $S\simeq\operatorname{Inn}(S)\leq G\leq
\operatorname{Aut}(S)$ for a nonabelian finite simple group $S$). The proof of Theorem~1  uses the following
statement on the number of classes of $\pi$-Hall subgroups in finite simple groups. 

\begin{Theo} \label{Simple} {\em ([3, Theorem~1.1], mod CFSG)}
Let $\pi$ be a set of primes and $G$ be an almost simple
finite $E_\pi$-group with the {\rm (}nonabelian simple{\rm )} socle~$S$. Then the following hold.

$(1)$~If $2\not\in\pi$ then  $k_\pi^G(S)=1.$

$(2)$~If $3\not\in\pi$ then $k_\pi^G(S)\in\{1,2\}$.

$(3)$~If $2,3\in\pi$ then $k_\pi^G(S)\in\{1,2,3,4,9\}$.

In particular, $k_\pi^G(S)$ is a $\pi$-number..
\end{Theo}

\begin{Lemma} Let $H$ be a $\pi$-Hall subgroup, let $A$ be a normal subgroup of  $G$, and $HAC_G(A)$ is normal in $G$ (this
condition is satisfied, if $HA\unlhd G$). Then an $A$-class of  $\pi$-Hall subgroups is a class of $G$-induced
$\pi$-Hall subgroups if and only if it is $H$-invariant. \end{Lemma}

\begin{proof} If $K$ is a $\pi$-Hall subgroup of $G$ then  $K\leq HAC_G(A)$, and so
$KAC_G(A)=HAC_G(A)$. Since the $A$-class  $\{(K\cap A)^a\mid a\in A\}$ is $K$-invariant, it follows that it is
invariant under  $HAC_G(A)=KAC_G(A)$, and so under~$H$.

Conversely, without loss of generality we may assume that $G=HA$ and the claim follows from Lemma~6. \end{proof}

\begin{Lemma} Let $H$ be a $\pi$-Hall subgroup, let $A$ be a normal subgroup of $G$, and $HA\unlhd G$. Then
$k_\pi^G(A)=k_\pi^{HA}(A)$. \end{Lemma}

\begin{proof} Since $HA$ is a normal subgroup of  $G$, each $\pi$-Hall subgroup of  $G$ lies in
$HA$. Hence~${k_\pi^G(A)=k_\pi^{HA}(A)}$. \end{proof}

\begin{Lemma} Let $H$ be a $\pi$-Hall subgroup, let $A$ be a normal subgroup of $G$, and $HA\unlhd G$. Then the
following  are equivalent:

$(1)$~$k_\pi^G(A)=1$.

$(2)$~$HA\in C_\pi$.

$(3)$~Every two $\pi$-Hall subgroups of $G$ are conjugate by an element of~$A$. \end{Lemma}

\begin{proof} $(1)\Rightarrow (2)$. If $K$ is a $\pi$-Hall subgroup of $HA$ then by (1) $H\cap A$ and
$K\cap A$ are conjugate in  $A$. We may assume that  $H\cap A=K\cap A$. Then   $H$ and $K$ lie in
$N_{HA}(H\cap A)$. By the Frattini argument,  $HA=N_{HA}(H\cap A)A$. So, $$N_{HA}(H\cap
A)/N_A(H\cap A)=N_{HA}(H\cap A)/N_{HA}(H\cap A)\cap A\simeq N_{HA}(H\cap A)A/A=HA/A$$ is a $\pi$-group. Thus
$N_{HA}(H\cap A)$ possesses a normal series $$N_{HA}(H\cap A)\geq N_A(H\cap A)\geq H\cap A\geq 1$$ such that
every factor of the series is either a $\pi$- or $\pi'$- group, and, by Lemma~4(2), satisfies $D_\pi$. In
particular, $H$ and $K$ are conjugate in~${N_{HA}(H\cap A)}$.

 $(2)\Rightarrow (3)$ and  $(3)\Rightarrow (1)$ are evident.
\end{proof}

\begin{Lemma} Let $H$ be a $\pi$-Hall subgroup, let $A=A_1\times\ldots\times A_s$ be a normal subgroup of $G$, and
$G=HAC_G(A)$. Then for every $i=1,\ldots,s$ the following  hold:

{\rm (1)}~$N_G(A_i)=N_H(A_i)AC_G(A)$.

{\rm (2)}~$N_H(A_i)$ is a $\pi$-Hall subgroup of~$N_G(A_i)$.

{\rm (3)}~$k_\pi^{\operatorname{Aut}_G(A_i)}(\operatorname{Inn}(A_i))=k_\pi^{N_G(A_i)}(A_i)$.
\end{Lemma}

\begin{proof}  (1) follows since  $G=HAC_G(A)$ and  $AC_G(A)\leq N_G(A_i)$. 
Using  (1) and the identity $N_H(A_i)\cap AC_G(A)=H\cap AC_G(A)$, we see that $\vert
N_G(A_i):N_H(A_i)\vert=\vert AC_G(A):(H\cap AC_G(A))\vert$ is a  $\pi'$-number, whence~(2).

Assume that $\rho:A_i\rightarrow \operatorname{Inn}(A_i)$ is the natural epimorphism. Since
$\operatorname{Ker}(\rho) =Z(A_i)$ is an abelian group, the kernel of $\rho$ possesses a unique
$\pi$-Hall subgroup which lies in each  $\pi$-Hall subgroup of~$A_i$. Hence the map $H\mapsto H\rho$
defines a bijection between the sets of $\pi$-Hall subgroups of  $A_i$ and $\operatorname{Inn}(A_i)$, and also
induces a bijection (we denote it by the same symbol $\sigma$) between the sets $\Delta$ and $\Gamma$ of $A_i$-
and $\operatorname{Inn}(A_i)$- classes of $\pi$-Hall subgroups, respectively. We show that the restriction of
$\sigma$ on the set $\Delta_0$ of all  $A_i$-classes of $N_G(A_i)$-induced $\pi$-Hall subgroups is a bijective
map from $\Delta_0$ onto the set $\Gamma_0$ of all $\operatorname{Inn}(A_i)$-classes of
$\operatorname{Aut}_G(A_i)$-induced $\pi$-Hall subgroups. Since  $A=A_1\times\ldots\times A_s$; therefore, (1)
implies $N_G(A_i)=N_H(A_i)A_iC_G(A_i)$. The normalizer $N_G(A_i)$ permutes elements from $\Delta$
acting by conjugation on the $\pi$-Hall subgroups of $A_i$. Thus, it acts on  $\Delta$. By Lemma 11, $\Delta_0$ is
the union of all one-element orbits under this action. By using  $\sigma$, define an equivalent action of
$N_G(A_i)$ on $\Gamma$. Since $C_G(A_i)$ lies in the kernel of both actions, the induced actions of
$\operatorname{Aut}_G(A_i)=N_G(A_i)/C_G(A_i)$ on $\Delta$ and $\Gamma$ are well-defined. It is easy to see that
the action of $\operatorname{Aut}_G(A_i)$ on $\Gamma$ defined in this way coincides with the natural action of
the group on the set of $\operatorname{Inn}(A_i)$-classes of conjugate $\pi$-Hall subgroups. Since
$$\operatorname{Aut}_G(A_i)/\operatorname{Inn}(A_i)\simeq N_G(A_i)/A_iC_G(A_i)\simeq N_H(A_i)/(N_H(A_i)\cap A_i
C_G(A_i))$$ is a $\pi$-group, by Lemma 11, $\Gamma_0$ coincides with the union of one-element orbits of
$\operatorname{Aut}_G(A_i)$ on $\Gamma$. By the definition of the action, $\Gamma_0$ is the image of $\Delta_0$
under $\sigma$. Since $\sigma$ is a bijection,
$$k_\pi^{N_G(A_i)}(A_i)=\vert\Delta_0\vert=\vert\Gamma_0\vert=k_\pi^{\operatorname{Aut}_G(A_i)}(\operatorname{Inn}(A_i)).$$
(3) follows.\end{proof}

Suppose $A=A_1\times\dots\times A_s$ and for every $i=1,\dots,s$ by ${\mathcal K}_i$ we denote an
$A_i$-class of $\pi$-Hall subgroups of $A_i$. The set  $$ {\mathcal K}_1\times \dots\times {\mathcal K}_s=\{\langle
H_1,\dots, H_s\rangle\mid H_i\in {\mathcal K}_i, \ i=1,\dots,s\} $$ is called the {\it product of classes} ${\mathcal
K}_1, \dots, {\mathcal K}_s$. Clearly ${\mathcal K}_1\times \dots\times {\mathcal K}_s$ is an $A$-class of $\pi$-Hall
subgroups of $A$. It is also clear that for a normal subgroup $A$ of  $G$ every $A$-class of $G$-induced
$\pi$-Hall subgroups is a product of some $A_1$-, $\dots$, $A_s$- classes of $G$-induced $\pi$-Hall subgroups. In
particular, $k_\pi^G(A_i)\le k_\pi^G(A)$ for every $i=1,\ldots,s$. The reverse inequality fails in general. 

\begin{Lemma} Let $H$ be a $\pi$-Hal subgroupl and let $A=A_1\times\dots\times A_s$ be a normal subgroup of $G$.
Assume also that the subgroups $A_1,\ldots,A_s$ are normal in $G$ and $G=HA C_G(A)$. Then
$k_\pi^G(A)=k_\pi^G(A_1)\cdot \ldots\cdot k_\pi^G(A_s)$. \end{Lemma}

\begin{proof} Two $\pi$-Hall subgroups $P$ and $Q$ of $A$ are conjugate in $A$ if and only if  $\pi$-Hall
subgroups $P\cap A_i$ and $Q\cap A_i$ of $A_i$  are conjugate in  $A_i$ for every $i=1,\dots,s$. In order to
prove the claim it is enough to show that the product of  $A_1$-, $\dots$, $A_s$- classes of $G$-induced
$\pi$-Hall subgroups is  an  $A$-class of  $G$-induced  $\pi$-Hall subgroups as well. Assume that $U_1,\dots,U_s$ are
$G$-induced $\pi$-Hall subgroups of $A_1,\dots,A_s$, respectively. We show that $$U=\langle
U_1,\dots,U_s\rangle=U_1\times\ldots\times U_s$$ is a $G$-induced $\pi$-Hall subgroup of $A$. By Lemma~11  it is
enough to show that for every $h\in H$ there exists $a\in A$ with $U^h=U^a$. Since $U_i=K_i\cap A_i$ for
appropriate $\pi$-Hall subgroup $K_i$ of $G$, the set $\bigl\{U_i^{x_i}\mid x_i\in A_i\bigr\}$ is
$K_i$-invariant, so this set is invariant under $K_iA=HA$. In particular,  $U_i^h=U_i^{a_i}$ for some $a_i\in
A_i$. Thus $$ U^h=U_1^h\times\ldots\times U_s^h= U_1^{a_1}\times\ldots\times U_s^{a_s}=U_1^a\times\ldots\times
U_s^a=U^a, $$ where $a=a_1\ldots a_s\in A$. \end{proof}

\begin{Lemma} Let $H$ be a $\pi$-Hall subgroup, let $A=A_1\times\dots\times A_s$ be a normal subgroup of $G$ that acts
transitively on the set $\{A_1,\dots, A_s\}$ by conjugation, and  $G=HA C_G(A)$. Then
$k_\pi^G(A)=k_\pi^G(A_i)=k_\pi^{N_G(A_i)}(A_i)$ for every~${i=1,\ldots,s}$. \end{Lemma}

\begin{proof} We show that each $G$-induced $\pi$-Hall subgroup of $A_i$ is also a $N_G(A_i)$-induced subgroup.
Indeed, if $K$ is a  $\pi$-Hall subgroup of  $G$, then $G=KAC_G(A)$ and, by Lemma 14, the identity
$N_G(A_i)=N_K(A_i)AC_G(A)$ holds. Moreover $K\cap A_i=N_K(A_i)\cap A_i$ and $N_K(A_i)$ is a $\pi$-Hall subgroup
of $N_G(A_i)$. So every $G$-induced $\pi$-Hall subgroup of $A_i$ is also a $N_G(A_i)$-induced $\pi$-Hall
subgroup, in particular, $k_\pi^G(A_i)\le k_\pi^{N_G(A_i)}(A_i)$.

Now we show that if $x,y\in G$ are in the same coset of  $G$ by $N_G(A_1)$, then for every $G$-induced $\pi$-Hall
subgroup $U_1$ of $A_1$ the subgroups $U_1^x$ and $U_1^y$ are conjugate in $A_i=A_1^x=A_1^y$. It is enough to
show that the subgroups $U_1$ and $U_1^t$, where $t=xy^{-1}\in N_G(A_1)$, are conjugate in $A_1$. Put $t=hac$,
$a\in A$, $c\in C_G(A)$, $h\in H$. Since  $U_1^h$ and $U_1^{hac}$ are conjugate in $A_1$, it is
enough to show that $U_1$ and $U_1^h$ are conjugate in  $A_1$. Since $(ac)^{h^{-1}}\in AC_G(A)\leq N_G(A_1)$, the
element $h$ normalizes $A_1$ as well. Assume that $U_1=U\cap A_1$ for a $G$-induced $\pi$-Hall subgroup $U$  of $A$.
Suppose that ${\mathcal K}$ is an $A$-class of $\pi$-Hall subgroups containing $U$, and let ${\mathcal K}={\mathcal K}_1\times
\dots\times {\mathcal K}_s$, where ${\mathcal K}_i$ is an $A_i$-class of $G$-induced $\pi$-Hall subgroups. Clearly,
$U_1\in {\mathcal K}_1$. Since, by Lemma  11,  ${\mathcal K}$ is $H$-invariant,  $H$ acts on the set
$\{{\mathcal K}_1, \dots, {\mathcal K}_s\}$. The element $h$ normalizes $A_1$, hence it fix the $A_1$-class ${\mathcal K}_1$.
In particular, the subgroups $U_1$ and $U_1^h$ are in ${\mathcal K}_1$, and so they are conjugate in~$A_1$.

Assume that $f\in H$ and $A_1^f=A_i$. For an $A_1$-class of $\pi$-Hall subgroups ${\mathcal K}_1$ we define the
$A_i$-class ${\mathcal K}_1^{f}$ by ${\mathcal K}_1^{f}=\{U_1^{f}\mid U_1\in {\mathcal K}_1\}$. As we noted above, ${\mathcal
K}_1^{f}$ is an $A_i$-class of $G$-induced $\pi$-Hall subgroups.

Let $h_1=1,h_2\ldots, h_s$ be the right transversal of $N_H(A_1)$ in $H$. Since $G$ acts transitively, up to
renumbering we may assume that $A_i=A_1^{h_i}$ and  $(N_G(A_1))^{h_i}=N_G(A_i)$. So $k_\pi^{N_G(A_1)}(A_1)=
k_\pi^{N_G(A_i)}(A_i)$ for $i=1,\ldots,s$. Consider  $$\sigma:{\mathcal K}_1\mapsto {\mathcal
K}_1^{h_1}\times\dots\times{\mathcal K}_1^{h_s},$$ mapping an $A_1$-class of $N_G(A_1)$-induced $\pi$-Hall subgroups
$\mathcal{K}_1$ to an $A$-class of $\pi$-Hall subgroups. Note that $\mathcal{K}_1^{h_1}\times\ldots\times
\mathcal{K}_1^{h_s}$ is always $H$-invariant and, by Lemma 11, it is an $A$-class of $G$-induced $\pi$-Hall
subgroups. Notice also that  $\sigma$ is injective, and so $k_\pi^{N_G(A_1)}(A_1)\le k_\pi^G(A)$. Consider the restriction $\tau$ of
$\sigma$ on the set of $A_1$-classes of $G$-induced $\pi$-Hall
subgroups. We need to show that the image of $\tau$ coincides with the set of $A$-classes of $G$-induced
$\pi$-Hall subgroups in order to complete the prove, since in such case we derive the
inequality~${k_\pi^G(A_1)\ge k_\pi^G(A)}$.

Let ${\mathcal K}={\mathcal K}_1\times \dots\times {\mathcal K}_s$ be an $A$-class of $G$-induced $\pi$-Hall subgroups. It is
enough to show that ${\mathcal K}_i={\mathcal K}_1^{h_i}$ for every $i=1,\dots, s$. Since $G$ acts transitively on the
set $\{A_1,\dots,A_s\}$, there exists an element $g\in G$ such that $A_1^g=A_i$. Let $gct$, where
$a\in A$, $c\in C_G(A)$, and $t\in H$. Then $A_1^t=A_i$ and $t\in N_H(A_1)h_i$. As we already proved, ${\mathcal
K}_1^t={\mathcal K}_1^{h_i}$. By Lemma 11, ${\mathcal K}$ is $H$-invariant. Hence ${\mathcal K}_i={\mathcal K}_1^t={\mathcal
K}_1^{h_i}$ and ${\mathcal K}={\mathcal K}_1^{h_1}\times\dots\times{\mathcal K}_1^{h_s}={\mathcal K}_1\tau.$\end{proof}

\section{A conjugacy criterion for Hall subgroups}

In the section we prove Theorem~1, Corollary~2, and provide an algorithm for determining whether $G$ satisfies
$C_\pi$ by using a normal series of~$G$.\medskip 

\noindent{\itshape Proof of Theorem 1.} Assume that the claim is not true, and let $G$ be a counterexample of minimal order.
Then  $G$ possesses a $\pi$-Hall subgroup $H$ and a normal subgroup $A$ such that $HA$ does not satisfy $C_\pi$.
We choose $A$ to be minimal. Let  $K$ be a $\pi$-Hall subgroup of $HA$ that is not conjugate with  $H$ in $HA$.
We divide into several steps the process of canceling the group~$G$.

Clearly

$(1)$ $HA=KA$.

$(2)$ $A$  {\sl is a minimal normal subgroup of $G$.}

Otherwise assume that $M$ is a nontrivial normal subgroup of $G$ that is contained properly in $A$. Put
$\overline{G}=G/M$ and, given a subgroup $B$ of $G$, denote  $BM/M$ by $\overline{B}$. By Lemma~9,
$\overline{G}$ satisfies $C_\pi$, $\overline{H}$ and $\overline{K}$ are  $\pi$-Hall subgroups of $\overline{G}$,
$\overline{A}$ is a normal subgroup of $\overline{G}$, $\overline{H}\overline{A}=\overline{K}\overline{A}$ and
$|\overline{G}|<|G|$. In view of the minimality of $G$, the group $\overline{H}\overline{A}$ satisfies $C_\pi$.
So  $\overline{H}$ and $\overline{K}$ are conjugate by an element of $\overline{A}$. Hence, the
subgroups  $HM$ and $KM$ are conjugate by an element of $A$. Without loss of generality, we may assume that
$HM=KM$. In view of the choice of $A$, the group $HM$ satisfies $C_\pi$. Hence $H$ and $K$ are conjugate by an
element of $M\leq A$; a contradiction.

$(3)$ $A\not\in C_\pi$. {\sl In particular, $A$ is not solvable.}

Otherwise, by Lemma~5, the group $HA$ satisfies $C_\pi$ as an extension of a $C_\pi$-groups by a $\pi$-group.

$(4)$ $HA$ {\sl is a normal subgroup of $G$.}

Otherwise  $N_G(HA)$ is a proper subgroup of $G$ and, by Lemma~7, we have $N_G(HA)\in C_\pi$. In view of the
minimality of $G$, it follows that  $HA\in C_\pi$; a contradiction.

 In view of $(2)$ and $(3)$

$(5)$
 $A$ {\sl is a direct product of simple nonabelian groups $S_1,\dots,S_m$. The group $G$ acts transitively on the
 set $\Omega=\{S_1,\dots,S_m\}$ by conjugation.}

Let $\Delta_1,\dots,\Delta_s$ be the orbits of $HA$ on $\Omega$, and put $T_j=\langle\Delta_j\rangle$ for
every $j=1,\dots,s$. In view of $(4)$ and $(5)$,

$(6)$ $G$ {\sl acts transitively on  $\{T_1,\dots, T_s\}$ by conjugation. The subgroup $A$ is a direct product of
$T_1,\dots,T_s$ and each of these subgroups is normal in~$HA$.}

Assume that $S\in\Omega$, and let $T$ be a subgroup, generated by the orbit from the set $\{\Delta_1,\dots,\Delta_s\}$
that contains $S$. By Lemma~16,

$(7)$ $k_\pi^{HA}(T)=k_\pi^{HA}(S)$.

By $(7)$ and Lemmas~12  and~15

$(8)$ $k_\pi^G(A)=k_\pi^{HA}(A)\bigl(k_\pi^{HA}(T)\bigr)^s=\bigl(k_\pi^{HA}(S)\bigr)^s$.

By $(8)$, Theorem~10 and Lemmas~14 and~16

$(9)$ $k_\pi^G(A)$  {\sl is a $\pi$-number.}

By Lemma~11,

$(10)$ $HA$ {\sl fixes every  $A$-class of $G$-induced $\pi$-Hall subgroups.}

Since  $G\in C_\pi$,

$(11)$ $G$ {\sl acts transitively on the set of $A$-classes of $G$-induced $\pi$-Hall subgroups.}

In view of $(10)$, the subgroup $HA$ lies in the kernel of this action. Now, by $(11)$,

$(12)$ $k_\pi^G(A)$ {\sl is a $\pi'$-number.}

By $(9)$ and $(12)$, it follows that

$(13)$ $k_\pi^G(A)=1$.

Now by Lemma~13,

$(14)$ $HA\in C_\pi$; {\sl a contradiction.} \qed \medskip

\noindent{\itshape Proof of Corollary 2. Necessity.} If $G\in C_\pi$, then, by Lemma 9, $G/A\in C_\pi$ as well. Let $K/A$
be a  $\pi$-Hall subgroup of $G/A$. By Lemma 8, there exists a $\pi$-Hall subgroup $H$ of $G$ such that $K=HA$.
By Theorem 1, the group $K=HA$ satisfies~$C_\pi$.

\noindent{\itshape Sufficiency.} Let $H$ be a $\pi$-Hall subgroup of $K$. Since  $K/A$ is a $\pi$-Hall subgroup of $G/A$, we
have $\vert G:H\vert=\vert G:K\vert\cdot\vert K:H\vert$ is a $\pi'$-number, and so $H$ is a $\pi$-Hall subgroup of
$G$. In particular, $G\in E_\pi$. Let $H_1,H_2$ be $\pi$-Hall subgroups of $G$. Since $G/A\in C_\pi$, the
subgroups $H_1A/A$ and $H_2A/A$ are conjugate in $G/A$, and we may assume that  $H_1A=H_2A$. However $H_1A\in C_\pi$,
and so $H_1$ and $H_2$ are conjugate. Thus ${G\in C_\pi}$. \qed\medskip

\begin{Lemma} Assume that $G=HA$, where $H$ is a $\pi$-Hall subgroup,  $A$ is a normal subgroup of $G$, and
$A=S_1\times\ldots\times S_k$ is a direct product of simple groups. Then $G\in C_\pi$ if and only if
$\Aut_G(S_i)\in C_\pi$ for every $i=1,\ldots,k$. \end{Lemma}

\begin{proof} By the Hall theorem and Lemma 5 we may assume that  $S_1,\ldots,S_k$ are nonabelian simple groups,
and so the set $\{S_1,\ldots,S_k\}$ is invariant under the action of $G$ by conjugation. Moreover this set is
partitioned into the orbits $\Omega_1,\ldots,\Omega_m$. Denote an element of  $\Omega_j$ by  $S_{i_j}$. By Lemmas
15 and 16, we obtain

\begin{equation} k_\pi^G(A)=k_\pi^G(S_{i_1})\cdot\ldots\cdot k_\pi^G(S_{i_m}).
\end{equation}

Lemmas 12 and 13 imply that $HA\in C_\pi$ if and only if $k_\pi^G(A)=1$, and in view of  (1), if and only if
$k_\pi^G(S_i)=1$ for every $i$. By Lemmas 14 and 16,
$k_\pi^G(S_i)=k_\pi^{N_G(S_i)}(S_i)= k_\pi^{\Aut_G(S_i)}(S_i)$ for every $i$. Moreover, since $\vert
\Aut_G(S_i):S_i\vert$ is a $\pi$-number, in view of Lemma 13,  $k_\pi^{\Aut_G(S_i)}(S_i)=1$ holds if
and only if $\Aut_G(S_i)\in C_\pi$. Therefore,  $k_\pi^G(S_i)=1$ for every $i$  if and only
if~${\Aut_G(S_i)\in C_\pi}$. \end{proof}

Now we are able to give an algorithm reducing the problem, whether a finite group satisfies $C_\pi$, to the check
of $C_\pi$-property in some almost simple groups. Assume  that 
\begin{equation}
 G=G_0>G_1>\ldots>G_n=1
 \end{equation}
is a chief series of $G$. Put $H_1=G=G_0$. Suppose that for some $i=1,\ldots,n$ the group  $H_i$ is
constructed so that $G_{i-1}\leq H_i$ and $H_i/G_{i-1}$ is a $\pi$-Hall subgroup of $G/G_{i-1}$. Since  (2) is a
chief series,   $$ G_{i-1}/G_i=S_1^i\times\ldots\times S_{k_i}^i, $$  where $S_1^i,\ldots,
S_{k_i}^i$ are simple groups. We check whether $$ \Aut_{H_{i}}(S^i_1)\in C_\pi, \ldots,\Aut_{H_{i}}(S^i_{k_i})\in
C_\pi. $$ If so then by Lemma~17, we have $H_{i}/G_{i}\in C_\pi$ and we can take a complete preimage of a
$\pi$-Hall subgroup of $H_{i}/G_i$ to be equal to $H_{i+1}$. Otherwise, by Corollary 2, we deduce that $G\not\in
C_\pi$ and stop the process. By Corollary 2 it follows, that $G$ satisfies $C_\pi$ if and only if the group
$H_{n+1}$ can be constructed. Notice that in this case $H_{n+1}$ is a $\pi$-Hall subgroup of~$G$. 

\begin{Cor} If either $2\not\in\pi$ or $3\not\in\pi$, then $G\in C_\pi$ if and only if every
nonabelian composition factor of $G$  satisfies~$C_\pi$. \end{Cor}

\begin{proof} The {\em sufficiency} follows from Lemma 5. We prove the {\em necessity}. By the above algorithm, we
may assume that $S\leq G\leq \Aut(S)$ for a nonabelian finite simple group $S$, and $G/S$ is a $\pi$-number. We
need to show that $S\in C_\pi$. Assume the contrary. Then Lemmas 11 and 13 imply that $G$ stabilizes precisely
one class of  $\pi$-Hall subgroups of $S$. Therefore, $S$ possesses at least three classes of $\pi$-Hall
subgroups. On the other hand, since either $2\not\in\pi$, or $3\not\in\pi$, by Theorem 10 the number of classes
of $\pi$-Hall subgroups in $S$ is not greater than $2$; a contradiction. \end{proof}

Notice that in the case  $2\not\in\pi$ the corollary is immediate from~[8, Theorem~A] and Lemma~5.

\end{document}